\documentclass[12pt]{article}
\usepackage{amscd}
\usepackage{amsthm}
\usepackage{extarrows}
\usepackage{latexsym}
\usepackage{mathrsfs}
\usepackage{pst-all,multido,ifthen}
\usepackage{amsmath}
\usepackage{amssymb}
\usepackage[all]{xy}
\usepackage[english]{babel}
\usepackage{mathtools}
\usepackage{upgreek}

\newtheorem{thm}{Theorem}
\newtheorem{df}{Definition}

\newtheorem{lemma}{Lemma}

\newtheorem{fact}{Fact}

\newtheorem{ex}{Example}
\def\endproof{$\hfill \square$}
\def\ad{\textup{ad}}
\def\Ad{\textup{Ad}}

\def\Aut{\textup{Aut}}
\def\char{\textup{char}}
\def\codim{\textup{codim}}
\def\End{\textup{End}}
\def\Ext{\textup{Ext}}
\def\Frac{\textup{Frac}}
\def\Fr{\textup{Fr}}
\def\Gal{\textup{Gal}}
\def\limind{\textup{lim.ind.}}
\def\GL{\textup{GL}}
\def\h{\textup{h}}
\def\sc{\textup{sc}}
\def\Hom{\textup{Hom}}
\def\GL{\textup{GL}}
\def\Im{\textup{Im}}
\def\Lat{\textup{Lat}}
\def\Lie{\textup{Lie}}
\def\Ker{\textup{Ker}}
\def\op{\textup{op}}
\def\PGL{\textup{PGL}}
\def\Rep{\textup{Rep}}
\def\SL{\textup{SL}}
\def\Spec{\textup{Spec}\,}

\usepackage{scalerel}
\usepackage{stackengine,wasysym}
\stackMath
\usepackage{tikz}
\newcommand*\circled[1]{\tikz[baseline=(char.base)]{
\node[shape=circle,draw,inner sep=0.1pt] (char) {#1};}}

\makeatletter
\newcommand{\medoplus}{\mathbin{\mathpalette\make@small\oplus}}
\newcommand{\medotimes}{\mathbin{\mathpalette\make@small\otimes}}

\newcommand{\make@small}[2]{%
 \vcenter{\hbox{%
 \scalebox{1.6}{$\m@th#1#2$}%
 }}%
}
\makeatother 
 
\begin{document}
\title{On Lie algebra modules which are modules over semisimple group schemes}
\author{Micah Loverro and Adrian Vasiu}
\maketitle

\noindent
{\bf Abstract.} Let $p$ be a prime. Given a split semisimple group scheme $G$ over a normal integral domain $R$ which is a faithfully flat $\mathbb Z_{(p)}$-algebra, we classify all finite dimensional representations $V$ of the fiber $G_K$ of $G$ over $K:=\Frac(R)$ with the property that the set of lattices of $V$ with respect to $R$ which are $G$-modules is as well the set of lattices of $V$ with respect to $R$ which are $\Lie(G)$-modules. We apply this classification to get a general criterion of extensions of homomorphisms between reductive group schemes over $\Spec K$ to homomorphisms between reductive group schemes over $\Spec R$. We also show that for a simply connected semisimple group scheme over a reduced $\mathbb Q$--algebra, the category of its representations is equivalent to the category of representations of its Lie algebra.



\section{Introduction}\label{S1}

Let $R$ be a commutative ring with $1\neq 0$. 

Let $r\in\mathbb N$. Let $G$ be a {\it semisimple group scheme}  over $\Spec R$ of rank $r$: it is an affine smooth group scheme over $\Spec R$ whose geometric fibers are semisimple groups over algebraically closed fields which admit maximal tori of dimension $r$. We recall that $G$ is called {\it split} if it has a maximal torus isomorphic to $\mathbb G_{m,R}^r$. If $R$ is connected, let $d\in\mathbb N$ be the relative dimension of $G$ over $\Spec R$.

By a $G$-module we mean a finitely generated projective $R$-module $M$ endowed with a homomorphism $\rho_M:G\rightarrow\pmb{\Aut}_M$, where $\pmb{\Aut}_M$ is the affine smooth group scheme over $\Spec R$ of linear automorphisms of $M$: if $S$ is an $R$-algebra, then
$$\pmb{\Aut}_M(S):=\{f:S\otimes_R M\rightarrow S\otimes_R M|f\;\textup{is a bijective}\;S\textup{-linear map}\}.$$ 
Thus, if $M=R^n$ for some $n\in\mathbb N\cup\{0\}$, then $\pmb{\Aut}_M=\pmb{\GL}_{n,R}$ is a general linear group scheme over $\Spec R$. If $P$ is another $G$-module, then by a $G$-module map between $M$ and $P$ we mean an $R$-linear map $f:M\rightarrow P$ such that for each $R$-algebra $S$  (equivalently, for each smooth $R$-algebra $S$) and every $g\in G(S)$, we have an identity 
\begin{equation}\label{EQ01}
1_S\otimes f\circ\rho_M(S)(g)=\rho_{P}(S)(g)\circ 1_S\otimes f:S\otimes_R M\rightarrow S\otimes_R P.
\end{equation} 
Let $\Rep(G)$ be the category of $G$-modules. 

Let $\mathfrak g:=\Lie(G)$ be the {\it Lie algebra} of $G$: it is a Lie algebra over $R$ which as an $R$-module is locally of finite rank (of rank $d$ if $R$ is connected). As $R$-modules, we identify $\mathfrak g=\Ker(G(R[x]/(x^2)\rightarrow G(R))$, where the $R$-epimorphism $R[x]/(x^2)\rightarrow R$ maps $x+(x^2)$ to $0$. We recall that the Lie bracket on $\mathfrak g=\Ker(G(R[x]/(x^2))\rightarrow G(R))$ is defined by taking the (total) differential of the commutator morphism $[,]:G\times_{\Spec R} G\rightarrow G$ at the identity section $\Spec R\to G$. 

By a $\mathfrak g$-module we mean a finitely generated  projective $R$-module $L$ equipped with a Lie algebra homomorphism $\varrho_L:\mathfrak g\rightarrow\mathfrak g\mathfrak l_R(L)$, where 
$$\mathfrak g\mathfrak l_R(L):=\{e:L\rightarrow L|e\;\textup{is an}\;R\textup{-linear map}\}$$ 
is equipped with the usual Lie bracket $[,]$: if $e_1,e_2\in\mathfrak g\mathfrak l_R(L)$, then we have $[e_1,e_2]:=e_1\circ e_2-e_2\circ e_1$. Thus $\mathfrak g\mathfrak l_R(L)$ is the Lie algebra over $R$ which is associated to the $R$-algebra $\End_R(L)$ and is identified with $\Lie(\pmb{\Aut}_L)$. If $J$ is another $\mathfrak g$-module, then by a $\mathfrak g$-module map between $L$ and $J$ we mean an $R$-linear map $f:L\rightarrow J$ such that for all $a\in\mathfrak g$ we have an identity $f\circ \varrho_L(a)=\varrho_{J}(a)\circ f:L\rightarrow J$. Let $\Rep(\mathfrak g)$ be the category of $\mathfrak g$-modules. 

We have a natural functor
\begin{equation}\label{EQ02}
Lie=Lie_G:\Rep(G)\rightarrow\Rep(\mathfrak g)
\end{equation}
that maps a $G$-module $M$ defined by the representation $\rho_M:G\rightarrow\pmb{\Aut}_M$ to 
$$\varrho_M:=\Lie(\rho_M):\mathfrak g=\Lie(G)\rightarrow\Lie(\pmb{\Aut}_M)=\mathfrak g\mathfrak l_M.$$
Note the typesetting difference: $\Lie(G)=\mathfrak g$ is a Lie algebra over $R$, $Lie=Lie_G$ is a functor, and $Lie(M)=Lie_G(M)$ is a representation of $\Lie(G)=\mathfrak g$.

If $R$ is a field and $M$ is a simple $G$-module, then the irreducible representation $\rho_M:G\rightarrow \pmb{\Aut}_M$ is called {\it infinitesimally irreducible} if the $\mathfrak g$-module $M$ is simple as well, see \cite{Bor1}, Sect. 6; one also calls $M$ an infinitesimally simple $G$-module.

One would like first to classify all the $\mathfrak g$-modules which are  $G$-modules, i.e., are isomorphic to objects in the image of the functor (\ref{EQ02}), and second to apply such a classification to obtain extension results from $\Spec K$ to $\Spec R$ for homomorphisms between reductive group schemes that are in line with the extension results obtained in \cite{V1}, Subsect. 4.3, \cite{V2} and \cite{V3}. 

Let $K:=\mathcal N_R^{-1}R$ be the total quotient ring of $R$, where $\mathcal N_R$ is the multiplicative set of nonzero divisors of $R$. If $R$ is an integral domain, then $K$ is a field and we will denote its characteristic by $\char(K)$ . 

Let $p$ be a prime. We are mainly interested in the following two situations:

\medskip
{\bf (i)} The ring $R$ is a $\mathbb Q$--algebra.

\smallskip
{\bf (ii)} The ring $R$ is a faithfully flat $\mathbb Z_{(p)}$-algebra (i.e., $K$ is a $\mathbb Q$--algebra and for each point $z\in \Spec R$, its residue field $k_z$ has characteristic either $0$ or $p$, and there exist such points $z$ with $\char(k_z)=p$). 

\medskip\noindent
In the situation (i) we have the following classical result which in essence is well-known (for instance, when $R=K$ is a field see \cite{Mi} and in the general case see \cite{SGA3-3}, Exp. XXIV, Prop. 7.3.1 which implies the surjectivity of the functor (\ref{EQ02}) on objects without the reduced assumption):

\begin{thm}\label{T1} 
We assume that $R$ is a reduced $\mathbb Q$--algebra and $G$ is simply connected. Then the functor (\ref{EQ02}) is an equivalence of categories.
\end{thm}

The goal of this paper is to obtain variants of Theorem \ref{T1} for the situation (ii). As Theorem \ref{T1} fails in the situation (ii) (see Theorem \ref{T2} below), one is led to consider a fixed nonzero $G_K$-module $V$ (so, if $R$ is an integral domain, $V$ is a finite dimensional $K$-vector space) and to study the natural map
\begin{equation}\label{EQ03}
Lie=Lie_G:\Lat_G(V)\rightarrow\Lat_{\mathfrak g}(V)
\end{equation}
induced by the functor $Lie$ and denoted in the same way, where $\Lat_G(V)$ (resp. $\Lat_{\mathfrak g}(V)$) is the set of {\it lattices} of $V$ with respect to $R$ which are $G$-modules (resp. are $\mathfrak g$-modules). Here and in what follows, by a lattice of $V$ with respect to $R$ we mean an $R$-submodule $L$ of $V$ which is a finitely generated projective $R$-module and for which the injective $K$-linear map $K\otimes_R L\rightarrow V$ is a bijection. If $L$ is a $G$-module, then $Lie(L)$ is $L$ but viewed as a $\mathfrak g$-module via the functor (\ref{EQ02}).

Let $G^{\sc}$ be the simply connected semisimple group scheme cover of $G$; so $V$ is also a $G^{\sc}_K$-module.

To study the map (\ref{EQ03}) we will assume that $R$ is an integral domain and that $\char(K)=0$. Let $\overline{K}$ be an algebraic closure of $K$. We recall from \cite{SGA3-3}, Exp. XXV, Thm. 1.1 that there exists a unique (up to ordering) product decomposition 
$$G^{\sc}_{\overline{K}}=\prod_{i=1}^n G_{i,\overline{K}}$$ 
such that each $G_{i,\overline{K}}$ has a simple adjoint group scheme $G^{\ad}_{i,\overline{K}}:=G_{i,\overline{K}}/Z(G_{i,\overline{K}})$ over $\Spec\overline{K}$, where $Z(G_{i,\overline{K}})$ is the center of $G_{i,\overline{K}}$; here $n\in\mathbb N$. 

For references to the standard facts recalled in this paragraph see Subsection \ref{S21}. As $\char(K)=0$, it is well-known that the $G^{\sc}_{\overline{K}}$-module $\overline{K}\otimes_K V$ is semisimple and hence we write it as a direct sum 
$$\overline{K}\otimes_KV=\oplus_{j=1}^m \overline{V}_j$$ 
of simple $G^{\sc}_{\overline{K}}$-modules; here $m\in\mathbb N$. Each $\overline{V}_j$ admits a tensor product decomposition 
$$\overline{V}_j=\otimes_{i=1}^n \overline{V}_{ij},$$ where every $\overline{V}_{ij}$ is a simple $G_{i,\overline{K}}$-module and where every element
$$(g_1,\ldots,g_n)\in G^{\sc}(\overline{K})=\prod_{i=1}^n G_{i,\overline{K}}(\overline{K})$$ 
acts on $\overline{V}_j$ in the usual tensorial way: for all $v_{1j}\in \overline{V}_{1j},\ldots,v_{nj}\in \overline{V}_{nj}$, it maps $v_{1j}\otimes v_{2j}\otimes\cdots\otimes v_{nj}$ to $g_1(v_{1j})\otimes g_2(v_{2j})\otimes\cdots\otimes g_n(v_{nj})$. Moreover, for all $i\in\{1,\ldots,n\}$, if a maximal torus $T_{i,\overline{K}}$ of a Borel subgroup $B_{i,\overline{K}}$ of $G^{\sc}_{i,\overline{K}}$ is given and if $r_i\in\mathbb N$  is the dimension of $T_{i,\overline{K}}$, then to $B_{i,\overline{K}}$ corresponds a basis $\omega_{i,1},\ldots, \omega_{i,r_i}$ of dominant weights of the group of characters 
$$X^*(T_{i,\overline{K}}):=\Hom(T_{i,\overline{K}},\mathbb G_{m,\overline{K}})\simeq \mathbb Z^{r_i}$$ 
and the representation $\overline{V}_{ij}$ is uniquely determined by its highest weight
$$w_{ij}=\sum_{l=1}^{r_i} c_{ijl}\omega_{i,l},$$
where each $c_{ijl}\in\mathbb Z_{\ge 0}$. We have $w_{ij}=0$, i.e., $c_{ij1}=\cdots=c_{ijr_i}=0$, if and only if $\overline{V}_{ij}$ is a trivial simple $G^{\sc}_{i,\overline{K}}$-module (equivalently, $\dim_{\overline{K}}(\overline{V}_{ij})=1$).

\begin{df}\label{df1}
We say that the nonzero $G_K$-module (or $G^{\sc}_K$-module) $V$ is $p$-latticed if for each $(i,j)\in\{1,\ldots,n\}\times\{1,\ldots,m\}$, for every $l\in\{1,\ldots, r_i\}$ we have $c_{ijl}\in\{0,\ldots,p-1\}$. 
\end{df}
 
\begin{thm}\label{T2}
We assume that $R$ is a normal integral domain which is a faithfully flat $\mathbb Z_{(p)}$-algebra. We consider the following two statements on the fixed nonzero $G_K$-module $V$:
 
\medskip\noindent
\circled{\textup{1}} The map $Lie:\Lat_G(V)\rightarrow\Lat_{\mathfrak g}(V)$ is a bijection.
 
\smallskip\noindent
\circled{\textup{2}} The $G_K$-module is $p$-latticed.

\medskip
Then the following three properties hold:

\medskip
{\bf (a)} The implication $\circled{\textup{2}}\Rightarrow \circled{\textup{1}}$ always holds.

\smallskip
{\bf (b)} We assume that there exists a discrete valuation ring $D$ of mixed characteristic $(0,p)$ which is a subring of $R$ such that $G$ is the pullback of a semisimple group scheme $G_D$ over $\Spec D$ and the $G_K$-module $V$ is the pullback of a $G_{\Frac(D)}$-module $V_{\Frac(D)}$, where $\Frac(D)=D[\frac{1}{p}]$ is the subfield of $K$ which is the field of fractions of $D$ (for instance, this holds if $G$ is split). Then the implication $\circled{\textup{1}}\Rightarrow \circled{\textup{2}}$ holds.

\smallskip
{\bf (c)} If $G$ is split, then we have an equivalence $\circled{\textup{1}}\Leftrightarrow \circled{\textup{2}}$.
\end{thm}

\begin{ex}\label{EX1} We assume that $G=\pmb{\SL}_{2,R}$ and $R$ is as in Theorem \ref{T2}. Then the map (\ref{EQ03}) is a bijection if and only if the $G_K$-module $V$ is a direct sum of simple $G_K$-modules of dimension at most $p$. 

For instance, suppose $R=\mathbb Z_{(p)}$ and $V$ is simple of dimension $p+1$, so it is the $p$-th symmetric power $V = \mathbb Qx^p \oplus \mathbb Qx^{p-1}y \oplus \cdots \oplus \mathbb Q y^p$ of the standard $G_{\mathbb Q}$-module $\mathbb Qx\oplus \mathbb Qy$ of rank 2 (here $x$ and $y$ are viewed as indeterminates). Then the map (\ref{EQ03}) is not surjective: consider the lattice 
$$L:=\mathbb Z_{(p)}x^p\oplus \mathbb Z_{(p)}x^{p-1}y\oplus\cdots\oplus \mathbb Z_{(p)}xy^{p-1}\oplus \frac{1}{p}\mathbb Z_{(p)}(x^p +  y^p)$$
of $V$ with respect to $\mathbb Z_{(p)}$. Let $T$ be the split torus of $G$ which normalizes both $\mathbb Z_{(p)}x$ and $\mathbb Z_{(p)}y$; it has rank $1$, i.e., $T\simeq \mathbb G_{m,\mathbb Z_{(p)}}$. The elements of the standard $\mathbb Z_{(p)}$-basis of $\mathfrak g$ map $(x,y)$ to $(y,0)$ or $(0,x)$ or $(x,-y)$ (respectively) and thus map $\frac{1}{p}(x^p+y^p)$ to elements of $L$. This implies that $L$ is a $\mathfrak g$-module. But $L$ is not a $T$-module and thus it is also not a $G$-module. 
\end{ex}

The highest weights of Definition \ref{df1} show up in the works of Curtis and Borel (see \cite{Bor1}, Sects. 6 and 7; see also \cite{C1} and \cite{C2} for original results under certain restrictions such as $p\ge 7$): they are precisely all the highest weights which in characteristic $p$ define infinitesimally irreducible representations (see \cite{Bor1}, Thms. 6.4 and 7.5 (iii)). 

Theorem \ref{T1} is proved in Section \ref{S3} based on the review of Section \ref{S2} that recalls classical properties of roots and of closed subgroup schemes of semisimple group schemes over $\Spec R$. Theorem \ref{T2} is proved in Section \ref{S5} based on the proof of Theorem \ref{T1}, on the mentioned works of Curtis and Borel, and on the following general result proved in Section \ref{S4}. 

\begin{thm}\label{T3}
Let $H$ be a semisimple group over an algebraically closed field $\kappa$. Let $P$ be an $H$-module such that the $\Lie(H)$-module $P$ is semisimple of the same length as the $H$-module $P$. Then the $H$-module $P$ is itself semisimple. 
\end{thm}

The following example shows that the ``same length" assumption of Theorem \ref{T3} is always necessary in positive characteristic.

\begin{ex}\label{EX2}
We assume $\char(\kappa)=p$. Let $0\rightarrow L_1\rightarrow Q\rightarrow L_2\to 0$ be a nonsplit short exact sequence of $H$-modules with $L_1$ and $L_2$ simple and $\dim_{\kappa}(Q)>2$ (see \cite{J}, Part 2, Ch. 7, Sects. 7.1 and 7.2 for general examples). For an $H$-module $V$, we consider the $H$-module $V^{(p)}$ defined by the representation which is the composite of the surjective Frobenius homomorphism $H\rightarrow H^{(p)}$ and the pullback $H^{(p)}\rightarrow \pmb{\GL}^{(p)}_V=\pmb{\GL}_{V^{(p)}}$ via the Frobenius endomorphism $Frob$ of $\Spec\kappa$ of the representation defining $V$; here $V^{(p)}:=\kappa \otimes_{Frob,\kappa} V$. Let $P:=Q^{(p)}$. Then $P$ is a trivial (hence semisimple) $\Lie(H)$-module of length $\dim_{\kappa}(Q)>2$. But $P$ is an idecomposable $H$-module of length $2$; this is so as $L_1^{(p)}$ and $L_2^{(p)}$ are simple $H$-modules and the short exact sequence $0\to L_1^{(p)}\to P\to L_2^{(p)}\rightarrow 0$ of $H$-modules is nonsplit (see \cite{J}, Part 2, Ch. 10, Prop. 10.16 for the injectivity of the natural pullback map $\Ext^1_H(L_1,L_2)\rightarrow\Ext^1_H(L_1^{(p)},L_2^{(p)})$).  
\end{ex}

By combining Theorem \ref{T2} with \cite{V3}, in Section \ref{S6} we prove the following theorem which is an application used in \cite{V4} to simplify the arguments of \cite{V1}, Subsect. 4.3 on extending homomorphisms between reductive group schemes in contexts related to integral models of Shimura varieties of Hodge type. 

\begin{thm}\label{T4}
We assume that $R$ is a normal integral domain and a faithfully flat $\mathbb Z_{(p)}$-algebra. Let $\mathcal G_K$ be a simply connected semisimple group over $\Spec K$. Let $\mathcal V$ be a $\mathcal G_K$-module which is $p$-latticed and let $\mathcal H_K:=\Im(\mathcal G_K\rightarrow\pmb{\Aut}_{\mathcal V})$. Let $\mathcal M$ be a lattice of $\mathcal V$ with respect to $R$ such that there exists a perfect symmetric bilinear form $\mathcal B:\mathcal M\times\mathcal M\rightarrow R$ which, over $K$, it is fixed by $\mathcal H_K$ and whose restriction to $\Lie(\mathcal H_K)\cap\End_R(\mathcal M)$ is a unit of $R$ times the Killing form (thus we have $p>2$, see \cite{V3}, Prop. 3.5 (a)). Then the schematic closure of the image $\mathcal H_K$ in $\pmb{\Aut}_{\mathcal M}$ is a semisimple group scheme $\mathcal H$ over $\Spec R$ whose simply connected semisimple group scheme cover $\mathcal G$ extends $\mathcal G_K$ and has the same Lie algebra $\Lie(\mathcal H_K)\cap\mathfrak g\mathfrak l_{\mathcal M}$ as either $\mathcal H$ or its adjoint $\mathcal G^{\ad}$ (i.e., the isogenies $\mathcal G\rightarrow\mathcal H\rightarrow\mathcal H^{\ad}=\mathcal G^{\ad}$ are \'etale). 
\end{thm}

Theorem \ref{T2} (c) was first obtained by the first author in the case when $R$ is noetherian while he was a graduate student. 

\section{A review}\label{S2}

In this section we assume that $G=G^{\sc}$ is simply connected and split and that $R$ is connected.

\subsection{The split context}\label{S21}
In this subsection we assume that $K$ is a field with $\char(K)=0$. Thus a group scheme $\triangle$ of finite type over $\Spec K$ is smooth  (for Cartier's Theorem, for instance, see \cite{DG}, Ch. II, Sect. 6, Thm. of Subsect. 1.1). Moreover, $\triangle$ is a (not necessarily connected) reductive group if and only if $\triangle$ is linearly reductive, i.e., each $\triangle$-module is semisimple (completely reducible), see \cite{F}, p. 178. Similarly, if $\triangle$ is semisimple, then its Lie algebra $\Lie(\triangle)$ is semisimple (this can be easily checked over $\overline{K}$), and Weyl's complete reducibility theorem implies that each $\Lie(\triangle)$-module is semisimple (see \cite{Bou1}, Ch. I, Subsect. 6.2, Thm. 2). Thus the categories $\Rep(G_K)$ and $\Rep(K\otimes_R \mathfrak g)=\Rep(\Lie(G_K))$ are semisimple abelian categories. 

Let $T_K$ be a maximal torus of $G_K$ which is split. Let $B_K$ be a Borel subgroup of $G_K$ that contains $T_K$. The Lie algebra $\Lie(T_K)$ is a split Cartan subalgebra of $K\otimes_R \mathfrak g$ and thus $K\otimes_R \mathfrak g$ is also split. Moreover, $\Lie(B_K)$ is a Borel subalgebra of $K\otimes_R \mathfrak g$. The simple $G_K$-modules are classified by the dominant weights of $T_K$ with respect to $B_K$ (see \cite{J}, Part 2, Ch. 2, Cor. 2.7) and the simple $K\otimes_R \mathfrak g$-modules are classified by the dominant weights of $\Lie(T_K)$ with respect to $\Lie(B_K)$ (see \cite{Bou2}, Ch. VIII, Sect. 7, Cor. 2). 

\subsection{Roots}\label{S22}
For centers of semisimple group schemes see \cite{SGA3-3}, Exp. XXII, Cor. 4.1.7.

As $R$ is connected, from \cite{SGA3-3}, Exp. XXV, Thm. 1.1 we get that:

\medskip\noindent
$\bullet$ There exists a unique (up to isomorphism) simply connected split semisimple group $G_{\mathbb Z}$ over $\Spec\mathbb Z$ such that $G=\Spec R\times_{\Spec\mathbb Z} G_{\mathbb Z}$.

\smallskip\noindent
$\bullet$ There exists a unique direct sum decomposition $G_{\mathbb Z}=\prod_{i=1}^n G_{i,\mathbb Z}$, such that each $G_{i,\mathbb Z}$ has an adjoint group scheme $G^{\ad}_{i,\mathbb Z}:=G_{i,\mathbb Z}/Z(G_{i,\mathbb Z})$ whose geometric fibers are simple, where $Z(G_{i,\mathbb Z})$ is the center of $G_{i,\mathbb Z}$.

\medskip
Defining $G_i:=\Spec R\times_{\Spec\mathbb Z} G_{i,\mathbb Z}$, we get a product decomposition 
$$G=\prod_{i=1}^n G_i$$ 
over $\Spec R$ and a product decomposition $G_K=\prod_{i=1}^n G_{i,K}$ over $\Spec K$. 

For $i\in\{1,\ldots,n\}$ let $T_{i,\mathbb Z}$ be a (split) maximal torus of a Borel subgroup scheme $B_{i,\mathbb Z}$ of $G_{i,\mathbb Z}$, let $T_i:=\Spec R\times_{\Spec\mathbb Z} T_{i,\mathbb Z}$ and $B_i:=\Spec R\times_{\Spec\mathbb Z} B_{i,\mathbb Z}$. Therefore, $T:=\prod_{i=1}^n T_i$ is a maximal torus of the Borel subgroup scheme $B:=\prod_{i=1}^n B_i$ of $G$. Let $B^{\op}$ be the Borel subgroup scheme of $G$ which is the opposite of $B$ with respect to $T$.

We identify $\mathbb G_{m,\mathbb Z}=\Spec \mathbb Z[x,x^{-1}]$ and $\Lie(\mathbb G_{m,\mathbb Z})=\mathbb Z$ in such a way that $1\in\mathbb Z=\Lie(\mathbb G_{m,\mathbb Z})$ gets identified with the $\mathbb Z$-linear map 
$$\Omega^1_{\mathbb Z[x,x^{-1}]/\mathbb Z}=\mathbb Z\frac{dx}{x}\rightarrow \mathbb Z$$ that maps $\frac{dx}{x}$ to $1$. We also identify $\Lie(\mathbb G_{m,R})=R\otimes_{\mathbb Z} \Lie(\mathbb G_{m,\mathbb Z})=R$. 

We can assume that, if $K$ is a field, then the choices made in Section \ref{S1} and Subsection \ref{S21} are compatible with our notation, i.e., for each $i\in\{1,\ldots,n\}$ the maximal torus $T_{i,\overline{K}}$ and the Borel subgroup $B_{i,\overline{K}}$ of $G_{i,\overline{K}}$ are indeed the extensions to $\Spec \overline{K}$ of $T_i$ and $B_i$ (respectively) and $T_K=\Spec K\times_{\Spec Z} T_{\mathbb Z}$ and $B_K=\Spec K\times_{\Spec Z} B_{\mathbb Z}$. As such, we identify
$$\mathcal W_{G_i}:=X^*(T_{i,\overline{K}})=X^*(T_{i,K})=X^*(T_i):=\Hom(T_i,\mathbb G_{m,R})\simeq \mathbb Z^{r_i},$$
$$\mathcal W_G:=X^*(T_{\overline{K}})=X^*(T_K)=X^*(T)=\Hom(T,\mathbb G_{m,R})=\oplus_{i=1}^n \mathcal W_{G_i}$$
and we speak about the monoid of dominant weights
$$\mathcal W_G^{\ge 0}:=\bigoplus_{i=1}^n (\oplus_{j=1}^{r_i} \mathbb Z_{\ge 0}\omega_{ij})\subset \oplus_{i=1}^n \mathcal W_{G_i}=\mathcal W_G$$
of $G$ with respect to the split maximal torus $T$ of the Borel subgroup scheme $B$ of $G$. 

We recall that, if $K$ is a field, then the (integral) weights of $K\otimes_R \mathfrak g$ are the elements of 
$$\mathcal W_{\mathfrak g}:=\Hom_{\mathbb Z}(\Lie(T_{\mathbb Z}),\Lie(\mathbb G_{m,\mathbb Z}))=\Hom_{\mathbb Z}(\Lie(T_{\mathbb Z}),\mathbb Z)$$inside 
$$\Hom_K(\Lie(T_K),\Lie(\mathbb G_{m,K}))=\Hom_K(\Lie(T_K),K),$$
with the dominant weights $\mathcal W_{\mathfrak g}^{\ge 0}$ being those with respect to $B_K$. We have a bijection
$$\mathcal L:\mathcal W_G\rightarrow\mathcal W_{\mathfrak g}$$
given by the rule $w\rightarrow Lie_{T_K}(w)$ which induces a bijection $\mathcal L:\mathcal W_G^{\ge 0}\rightarrow\mathcal W_{\mathfrak g}^{\ge 0}$ denoted in the same way.

Let $G_{\mathbb Q}:=\Spec\mathbb Q\times_{\Spec\mathbb Z} G_{\mathbb Z}$. If $K$ is a field, as the abelian categories $\Rep(G_K)$ and $\Rep(G_{\mathbb Q})$ are semisimple, from the classification of simple $G_K$-modules  or $G_{\mathbb Q}$-modules in terms of dominant weights (see \cite{J}, Part II, Ch. 2, Cor. 2.7), we get that these simple modules are absolutely simple (see \cite{J}, Part II, Ch. 2, Cor. 2.9). Thus, if $K$ is a field, then the pullback functor 
$$K\otimes_{\mathbb Q} \Rep(G_{\mathbb Q})\rightarrow\Rep(G_K)$$ 
is an equivalence between $K$-linear semisimple abelian categories. In particular, for each $G_K$-module $V$, there exists a unique (up to isomorphism) $G_{\mathbb Q}$-module $V_{\mathbb Q}$ such that the $G_K$-modules $V$ and $K\otimes_{\mathbb Q} V_{\mathbb Q}$ are isomorphic. 

\subsection{Closed subgroups}\label{S23}
In this subsection we also assume that $\Spec R$ is reduced. Let $\Phi(G,T)$ be the root system of $G$ with respect to $T$ and let $\Phi^+(G,T)$ be the set of positive roots of $\Phi(G,T)$ with respect to $B$. We have a disjoint union 
$$\Phi(G,T)=\Phi^+(G,T)\sqcup -\Phi^+(G,T)$$ as well as direct sum decompositions of $R$-modules
$$\mathfrak g=\Lie(T)\oplus_{\alpha\in \Phi(G,T)} \mathfrak g_{\alpha}\;\;\;\textup{and}\;\;\;\Lie(B)=\Lie(T)\oplus_{\alpha\in \Phi^+(G,T)} \mathfrak g_{\alpha},$$
where each $\mathfrak g_{\alpha}$ is the weight space of the $T$-module $\mathfrak g$ corresponding to the weight $\alpha$, with $\mathfrak g$ being viewed as a $G$-module (hence also as a $T$-module) via the adjoint representation $\Ad:G\rightarrow\pmb{\Aut}_{\mathfrak g}$. For each $\alpha\in\Phi(G,T)$ there exists a unique $\mathbb G_{a,R}$ closed subgroup scheme $U_{\alpha}$ of $G$ which is normalized by $T$ and whose Lie algebra is $\mathfrak g_{\alpha}$ (see \cite{SGA3-3}, Exp. XII, Sect. 1, Thm. 1.1 or \cite{J}, Part II, Ch. 1, Sects. 1.1 and 1.2); as $R$ is a reduced $\mathbb Q$--algebra, $\Lie(U_{\alpha})=\mathfrak g_{\alpha}$ implies that $U_{\alpha}$ is normalized by $T$ as one can easily check based on \cite{Bor2}, Ch. II, Sect. 7, Subsect. 7.1. As $U_{\alpha}\simeq\mathbb G_{a,R}$, $\mathfrak g_{\alpha}$ is a free $R$-module of rank $1$.

We recall from \cite{SGA3-3}, Exp. XII, Sect. 4,  Prop. 4.1.2 that the product morphism
$$\imath:(\prod_{\alpha\in -\Phi^+(G,T)} U_{\alpha})\times_{\Spec R} T\times_{\Spec R} (\prod_{\alpha\in \Phi^+(G,T)} U_{\alpha})\rightarrow G$$
is an open embedding whose image $U:=\Im(\iota)$ does not depend on the orderings of the first and third factor of the source of $\imath$, being in fact equal to the image $B^{\op}B$ of the product morphism $B^{\op}\times_{\Spec R} B\rightarrow G$. In particular, $\iota$ induces an isomorphism
\begin{equation}\label{EQ04}
\jmath:(\prod_{\alpha\in\Phi^+(G,T)} U_{\alpha})\times_{\Spec R} T\times_{\Spec R} (\prod_{\alpha\in -\Phi^+(G,T)} U_{\alpha})\to U.
\end{equation} 

\section{Proof of Theorem \ref{T1}}\label{S3}

In this section we assume that $R$ is a reduced $\mathbb Q$--algebra and $G=G^{\sc}$ is simply connected. 

We write $R=\limind_{\lambda\in\Lambda} R_{\lambda}$ as an inductive limit of finitely generated $\mathbb Z$-subalgebras of $R$ where $\Lambda$ is the set of finite subsets of $R$ and $R_{\lambda}$ is the $\mathbb Z$-subalgebra of $R$ generated by the elements of $\lambda$. 

In this paragraph we recall the essentially well-known property that there exists $\lambda_0\in\Lambda$ such that $G$ is the pullback of a simply connected semisimple group scheme $G_{\lambda_0}$ over $\Spec R_{\lambda_0}$. As group objects in a category are defined in terms of commutative diagrams, from \cite{G2}, Thm. (8.8.2) we get that there exists $\lambda_0\in\Lambda$ such that $G$ is the pullback of an affine group scheme $G_{\lambda_0}$ over $\Spec R_{\lambda_0}$ of finite type. Based on \cite{SGA3-1}, Exp. $VI_B$, Subsect. 10.9 and Prop. 3.9 we can assume that there exists an open subgroup scheme $G^0_{\lambda_0}$ of $G_{\lambda_0}$ whose fibers over $\Spec R_{\lambda_0}$ are connected. From this and the affineness part of \cite{G2}, Thm. (8.10.5) we get that we can assume that $G^0_{\lambda_0}$ is affine and hence we can also assume that $G_{\lambda_0}=G^0_{\lambda_0}$. Based on \cite{G2}, Thm. (11.2.6) we can assume that $G_{\lambda_0}$ is flat over $\Spec R_{\lambda_0}$. From \cite{SGA3-3}, Exp. XIX, Thm. 2.5 we get that there exist a largest open subscheme $S_{\lambda_0}$ of $\Spec R_{\lambda_0}$ such that $S_{\lambda_0}\times_{\Spec R_{\lambda_0}} G_{\lambda_0}$ is a semisimple group scheme over $S_{\lambda_0}$. From this and the isomorphism part of \cite{G2}, Thm. (8.10.5) we get first that we can assume that $S_{\lambda_0}=\Spec R_{\lambda_0}$ and second that we can assume $G_{\lambda_0}$ is simply connected.

For $\lambda_0\subset\lambda\in\Lambda$, let $G_{\lambda}:=\Spec R_{\lambda}\times_{\Spec R_{\lambda_0}} G_{\lambda_0}$ and let $\mathfrak g_{\lambda}:=\Lie(G_{\lambda})$. 

As $\mathbb Z$ is a universally Japanese ring (see \cite{G1}, Cor. (7.7.4)), each finitely generated $\mathbb Z$-algebra which is an integral domain has a normalization which is a finitely generated $\mathbb Z$-algebra and hence noetherian. Thus, if $R$ is a normal integral domain, the normalization of each $R_{\lambda}$ is noetherian, and hence $R$ is the inductive limit of such normal noetherian integral domains.

From \cite{G2}, Thm. (8.5.2) we get that the categories $\Rep(G)$ and $\Rep(\mathfrak g)$ are the inductive limits of the categories $\Rep(G_{\lambda})$ and $\Rep(\mathfrak g_{\lambda})$ (respectively) indexed by $\lambda\in\Lambda$, $\lambda\supset\lambda_0$. Thus to prove Theorem \ref{T1}, by replacing $R$ with $R_{\lambda}$ for some $\lambda\in\Lambda$, $\lambda\supset\lambda_0$, we can assume that $R$ is a finitely generated $\mathbb Z$-algebra, hence noetherian; hence, as $R$ is reduced, $K$ is a finite product of fields. 

To check that the functor (\ref{EQ02}) is faithful, by replacing $R$ by a direct factor of $K$ which is a field, we can assume that $R=K$ is a field of characteristic zero and this case is well-known (for instance, using graphs, this follows from \cite{Bor2}, Ch. II, Sect. 7, Subsect. 7.1).

Thus to prove Theorem \ref{T1} it suffices to show that the functor (\ref{EQ02}) is surjective on objects and on morphisms. To check this, we can work locally in the \'etale topology of $\Spec R$ (cf. Equation (\ref{EQ01})) and hence we can also assume that $R$ is connected and that (see \cite{SGA3-3}, Exp. XIX, Prop. 6.1) $G$ has a maximal torus which is split. 

\subsection{Surjectivity of the functor (\ref{EQ02}) on objects}\label{S31}

Though the surjectivity of the functor (\ref{EQ02}) on objects follows from \cite{SGA3-3}, Exp. XXIV, Prop. 7.3.1 without even assuming that the $\mathbb Q$-algebra $R$ is reduced, several parts of the proof included here are used in the proof of Theorem \ref{T2}.

Let $1_{\star}$ be the identity automorphism of a (group) scheme $\star$. 

Let $L$ be a $\mathfrak g$-module. To check that there exists a $G$-module $M$ such that we have $M=L$ as $\mathfrak g$-modules, we consider four disjoint cases in order to include several proofs, including simpler ones in the easier cases such as when $R$ is a field or a discrete valuation ring or a normal integral domain.

\smallskip\noindent
{\bf Case 1: $R=K$ is a field.} We include two proofs in this case. 

The first proof, slightly sketched, is well-known and relies on the classification of simple $G$-modules and simple $\mathfrak g$-modules. As the abelian category $\Rep(\mathfrak g)$ is semisimple, we can assume that $L$ is a simple $\mathfrak g$-module. Let $w\in\mathcal W_G^{\ge 0}$ be such that $\mathcal L(w)\in \mathcal W_{\mathfrak g}^{\ge 0}$ is the dominant weight with the property that $L$, up to isomorphism, is the simple $\mathfrak g$-module of highest weight $\mathcal L(w)$. Let $M$ be the simple $G$-module of highest weight $w\in\mathcal G^{\ge 0}$. It is an easy exercise to check that $\Lie(M)$ and $L$ are isomorphic $\mathfrak g$-modules. 

For the second proof we consider a $G$-module $V$ such that the representation $\rho_V:G\rightarrow\pmb{\Aut}_V$ is faithful, to be a viewed as a closed embedding. Thus the $\mathfrak g$-module $N:=Lie(V)\oplus L=V\oplus L$ is such that the representation $\varrho_N:\mathfrak g\rightarrow\End(N)$ is also faithful, to be viewed as an inclusion. Let $\mathcal T(N):=\oplus_{a,b\ge 0} N^{\otimes a}\otimes_K (N^*)^{\otimes b}$, where $N^*:=\Hom_K(N,K)$ is the dual of $N$ and $\dag ^{\otimes c}$, with $c\in\mathbb Z_{\ge 0}$, is the tensor product over $K$ of $c$-copies of the $K$-vector space $\dag$ ($\dag^{\otimes 0}:=K$). From \cite{Se}, Ch. VI, Sect. 5, Thm. 5.2 we get that there exists a finite subset $\mathcal F\subset\mathcal T(N)$ with the property that $\mathfrak g$ is the Lie subalgebra of $\End_K(N)$ that annihilates every tensor of $\mathcal F$. We can assume that the projection $\pi\in\End_K(N)=N\otimes_K N^*$ of $N$ on $L$ along $V$ belongs to $\mathcal F$. Let $E$ be the subgroup of $\pmb{\Aut}_N$ which fixes each tensor of $\mathcal F$. Let $E^0$ be the connected component of the identity element of $E$. We have $\Lie(E)=\mathfrak g$, thus also $\Lie(E^0)=\mathfrak g$. As $E$ fixes $\pi$, both $L$ and $V$ are $E^0$-modules. In particular, we have a representation $\sigma:E^0\rightarrow  \pmb{\Aut}_V$ with the property that $Lie(\sigma)$ is injective (due to the identity $\Lie(E^0)=\mathfrak g$). This implies that $\sigma$ induces an \'etale isogeny $\sigma:E^0\rightarrow\Im(\sigma)$ and moreover we have $\Lie(\Im(\sigma))=\Lie(G)=\mathfrak g\subset\End_K(V)$. From this and \cite{Bor2}, Ch. II, Sect. 7, Subsect. 7.1 we get that $\Im(\sigma)=G$ and hence we have an \'etale isogeny $E^0\rightarrow G$. As $G$ is simply connected, we conclude that $E^0\rightarrow G$ is an isomorphism. As $L$ is an $E^0$-module, we conclude that it is as well a $G$-module in such a way that $Lie(L)$ is the $\mathfrak g$-module $L$.  

\smallskip\noindent
{\bf Case 2: $R$ is a discrete valuation ring.} Let $k$ be the residue field of $R$; we have $\char(k)=0$. 

Let $V$ be the  $G_K$-module such that $Lie(V)=K\otimes_R L$, see Case 1. We identify $V=K\otimes_R L$ as $K$-vector spaces and let $M$ be the lattice of $V$ with respect to $R$ which, under the mentioned identification, gets identified with $L$. It suffices to show that $M$ is a $G$-module. 

Based on \cite{SGA3-1}, Exp. VIB, Rm. 11.11.1 we get the existence of a $G$-module $P$ such that the representation $\rho_P:G\rightarrow\pmb{\Aut}_P$ is a closed embedding. From \cite{J}, Part I, Ch. X, Lem. of Sect. 10.4 we get that there exists a lattice $M'$ of $V$ which is a $G$-module. As $G$ is split and simply connected and as $R$ is connected, the existence of $P$ (respectively $M'$) also follows by pullback to $\Spec R$ from the mentioned references applied over $\mathbb Z_{(p)}$ to $G_{\mathbb Z_{(p)}}:=\Spec \mathbb Z_{(p)}\times_{\Spec\mathbb Z} G_{\mathbb Z}$ (respectively to $G_{\mathbb Z_{(p)}}$ and a $V_{\mathbb Q}$ as in Subsection \ref{S22}). Let $Q:=P\oplus M$ and $Q':=P\oplus M'$; we have $K\otimes_R Q=K\otimes_R Q'$.

The representation $\rho_{Q'}:G\rightarrow \pmb{\Aut}_{Q'}$ is a closed embedding and thus, with the notation of Subsection \ref{S23}, for each $\alpha\in\Phi(G,T)$, $U_{\alpha}$ is also a closed subgroup scheme of $\pmb{\Aut}_{Q'}$ whose Lie algebra is identified with $\mathfrak g_{\alpha}$ via the faithful representation $\varrho_{Q'}=Lie(\rho_{Q'}):\mathfrak g\rightarrow\End_R(Q')$.

For $\alpha\in\Phi(G,T)$ let $\mathbb V_{\alpha}$ be the vector group scheme over $\Spec R$ whose group of $S$-valued points is $S\otimes_R \mathfrak g_{\alpha}$ for each $R$-algebra $S$. As $R$ is a $\mathbb Q$--algebra and $Q$ and $Q'$ are $\mathfrak g$-modules, we have homomorphisms
$$\eta_{\alpha}:\mathbb V_{\alpha}\rightarrow \pmb{\Aut}_Q\;\;\;\textup{and}\;\;\;\eta'_{\alpha}:\mathbb V_{\alpha}\rightarrow \pmb{\Aut}_{Q'}$$ which for each $R$-algebra $S$ map $x\in \mathbb V_{\alpha}(S)=S\otimes_R \mathfrak g_{\alpha}$ to the sums $\sum_{q=0}^{\infty} \frac{\varrho_Q(x)^q}{q!}$ and $\sum_{q=0}^{\infty} \frac{\varrho_{Q'}(x)^q}{q!}$ (respectively). These sums coincide as elements of $\End_K(V)$ and are finite sums as each $x$ acts nilpotently on $S\otimes_R Q$ and $S\otimes_R Q'$. The image of $\eta_{\alpha,K}=\eta_{\alpha,K}'$ and $U_{\alpha,K}$ have the same Lie algebras and thus they coincide, see \cite{Bor2}, Ch. II, Sect. 7, Subsect. 7.1. This implies that $\eta_{\alpha}'$ factors through a homomorphism $\zeta_{\alpha}:\mathbb V_{\alpha}\rightarrow U_{\alpha}$ which induces an isomorphism at the level of Lie algebras, and hence is \'etale. As $\mathbb G_a$ over each field of characteristic $0$ has no finite nontrivial subgroup, we deduce that the fibers of $\zeta_{\alpha}$ are isomorphisms, based on which we easily see that $\zeta_{\alpha}$ itself is an isomorphism. 

We get a homomorphism $\eta_{\alpha}\circ\zeta_{\alpha}^{-1}:U_{\alpha}\rightarrow \pmb{\Aut}_Q$, hence $Q$ is a $U_{\alpha}$-module.

We fix an identification $T=\mathbb G_{m,R}^r$ and with respect to it we speak about the $\mathbb G_{m,R}$ factors of $T$ (there exist $r$ such factors). If $F=\mathbb G_{m,R}$ is such a factor of $T$, then we have a direct sum decomposition $K\otimes_R Q=\oplus_{q\in\mathbb Z} W_q$ such that $F_K$ acts on $W_q$ via the $q$-th power of the identity character of $F_K$. The standard generator $x$ of $\Lie(F)$ acts on $W_q$ as the multiplication by $q$. As $Q$ is a $\mathfrak g$-module, we have $x(Q)\subset Q$. As $x(Q)\subset Q$ and as for distinct integers $q_1,q_2$ which are eigenvalues of $x$ acting on the $K$-vector space $K\otimes_R Q$, the difference $q_1-q_2$ is invertible in $k$, it is an easy exercise to check that we have a direct sum decomposition $Q=\oplus_{q\in\mathbb Z} Q\cap W_q$. This implies that $Q$ is an $F$-module. The resulting homomorphism $\eta_F:F\rightarrow\pmb{\Aut}_Q$ is a closed embedding as this is so over $\Spec K$, see \cite{V3}, Lem. 2.3.2 (b) and (c). The images $\Im(\eta_F)$ indexed by such factors $F$ of $T$ commute as this is so over $\Spec K$ and hence we get a product homomorphism $\eta_T:T\rightarrow\pmb{\Aut}_Q$ which over $\Spec K$ is a closed embedding. Again from \cite{V3}, Lem. 2.3.2 (b) and (c) we get that $\eta_T$ is a closed embedding. In particular, $Q$ is a $T$-module. 

From the last two paragraphs we get a product morphism 
\begin{equation}\label{EQ05}
\eta:(\prod_{\alpha\in\Phi^+(G,T)} U_{\alpha})\times_{\Spec R} T\times_{\Spec R} (\prod_{\alpha\in -\Phi^+(G,T)} U_{\alpha})\rightarrow \pmb{\Aut}_Q
\end{equation}
which is the product of the $\eta_{\alpha}$'s and $\eta_T$ and which is compatible with the representation $\rho_{K\otimes_R Q}:G_K\rightarrow \pmb{\Aut}_{K\otimes_R Q}$, in the sense that $\rho_{K\otimes_R Q}$ restricted to $U_K$ is $\eta_K\circ\jmath_K^{-1}$. 

The union $U_+:=G_K\cup U$ is an open subscheme of $G$ whose complement $C:=G\setminus U_+$, when endowed with the reduced structure, is a reduced closed subscheme of $G_k$ of dimension less than $d=\dim(G_k)$. Thus, as we have $\dim(G)=d+\dim(R)=d+1$, we get that $\codim_G(C)\ge 2$. 

From the lats two paragraphs we get the existence of a morphism 
$$\rho_{Q,U_+}:U_+\rightarrow \pmb{\Aut}_Q$$ 
which extends both $\rho_{K\otimes_R Q}$ and the composite 
\begin{equation}\label{EQ06}
\rho_{Q,U}:=\eta\circ\jmath^{-1}:U=\Im(\iota)\to \pmb{\Aut}_Q.
\end{equation}
As $\codim_G(C)\ge 2$ and the scheme $\Spec R$ is normal noetherian, from \cite{BLR}, Ch. 4, Sect. 4.4, Thm. 1 we get that the morphism $\rho_{Q,U_+}$ extends to a morphism $\rho_Q:G\rightarrow \pmb{\Aut}_Q$ which, as it extends $\rho_{K\otimes_R Q}$, is a homomorphism. So $Q$ is a $G$-module.

The projection of $Q$ on $M$ along $P$ is fixed by $G$ (as it is fixed by $G_K$). This implies that $\rho_Q$ induces a homomorphism $\rho_M:G\rightarrow \pmb{\Aut}_M$ that extends $\rho_V$, hence $M$ is a $G$-module. 

\smallskip\noindent
{\bf Case 3: $R$ is normal but neither a field nor a discrete valuation ring.} Let $\mathcal D$ be the set of all local rings of $R$ which are discrete valuation rings; we recall (for instance, see \cite{Ma}, Thm. 11.5) that, as $R$ is noetherian, $\mathcal D$ is nonempty and in fact we have  $R=\cap_{O\in\mathcal D} O$.

Let $M:=L$. From Case 2 we get the existence of an open subscheme $Y$ of $\Spec R$ which contains all points of $\Spec R$ of codimension in $\Spec R$ at most $1$ (i.e., the closed subscheme $\Spec R\setminus Y$ has codimension in $\Spec R$ at least 2) and for which we have a homomorphism $\rho_{M,Y}:G_Y\rightarrow Y\times_{\Spec R} \pmb{\Aut}_M$ between reductive group schemes over $Y$ with the property that for each $O\in\mathcal D$, the $O\otimes_R \mathfrak g$-module $O\otimes_R M$ is $O\otimes_R L$. Considering the closed embedding 
$$(\rho_{M,Y},1_{G_Y}): G_Y\rightarrow (Y\times_{\Spec R} \pmb{\Aut}_M)\times_Y G_Y,$$ 
from \cite{V4}, Prop. 5.1 we get that it extends uniquely to a closed embedding homomorphism $(\rho_M,1_G): G\rightarrow \pmb{\Aut}_M\times_{\Spec R} G$. The resulting homomorphism $\rho_M:G\rightarrow \pmb{\Aut}_M$ endows $M$ with the structure of a $G$-module. The fact that the $\mathfrak g$-module structure on $M$ is the same one as the one given by $M=L$ follows from the fact that this is so over $O$ for one (hence all) $O\in\mathcal D$. 

\smallskip\noindent
{\bf Case 4: $R$ is not normal.} Let $V$ be the  $G_K$-module such that we have $Lie(V)=K\otimes_R L$, see Case 1 applied to the direct factors of $K$ which are fields. Let $G_{\mathbb Z_{(p)}}$ and $G_{\mathbb Q}$ be as in Subsection \ref{S22} and let $M$ be as in Case 2. As $R$ is connected and $G=G^{\sc}$ is split, there exists a $G_{\mathbb Q}$-module $V_{\mathbb Q}$ such that the $G_K$-module $V$ is isomorphic to $K\otimes_{\mathbb Q} V_{\mathbb Q}$ (see end of Subsection \ref{S22}). This implies that there exist $G$-modules $P$ and $M'$ such that $\rho_P:G\to \pmb{\Aut}_P$ is a closed embedding and $V=K\otimes_R M'$ (they are obtained, to be compared with Case 2, by pullback from $\Spec \mathbb Z_{(p)}$ to $\Spec R$). Based on this, as in Case 2 we argue the existence of a product morphism $\eta$ as in Equation (\ref{EQ05}), and hence we get a morphism $\rho_{Q,U}:U=\Im(\iota)\to\pmb{\Aut}_Q$ as in Equation (\ref{EQ06}). The pullback of $\rho_{Q,U}$ to $\Spec K$ coincides with the restriction to $U_K$ of the homomorphism $\rho_{K\otimes_R Q}:G_K\to\pmb{\Aut}_{K\otimes_R Q}$. 

The product morphism 
$$\Theta:U\times_{\Spec R} U\to G$$ 
is surjective and smooth, in particular it is a faithfully flat morphism between affine schemes. We will use affine faithfully flat descent with respect to $\Theta$ to show that the morphism $\rho_{Q,U}$ extends to a morphism $\rho_Q:G\to\pmb{\Aut}_Q$ that extends $\rho_{K\otimes_R Q}$. We consider the two projections 
$$\Pi_1,\Pi_2: (U\times_{\Spec R} U)\times_G (U\times_{\Spec R} U)\to U\times_{\Spec R} U$$ defined by $\Theta$. The two composite morphisms  
$$\rho_{Q,U}\circ\Pi_1,\rho_{Q,U}\circ\Pi_2:(U\times_{\Spec R} U)\times_G (U\times_{\Spec R} U)\to\pmb{\Aut}_Q$$
coincide as this is so after pullback to $\Spec K$. This implies the existence of $\rho_Q$ with the desired property. 

As $\rho_{K\otimes_R Q}$ is  homomorphism we get that $\rho_Q$ is a homomorphism. As in the last paragraph of Case 2 we argue that $M$ is a $G$-module. The fact that the $\mathfrak g$-module structure on $M$ is the same one as the one given by $M=L$ follows from the fact that this is so over $K$.

\subsection{Surjectivity of the functor (\ref{EQ02}) on morphisms}\label{S32}

Let $f:L\rightarrow J$ be a morphism of $\Rep(\mathfrak g)$. Based on Subsection \ref{S31}, we know that there exist $G$-modules $M$ and $P$ such that $Lie(M)=L$ and $Lie(P)=J$, i.e., we have $M=L$ and $P=J$ as $R$-modules but, in connection to $f:L\rightarrow J$, we view them as $\mathfrak g$-modules. We denote also by $f:M\rightarrow P$ the $R$-linear map defined by $f$ and the identifications $M=L$ and $P=J$, and to end the proof of Theorem \ref{T1} it suffices to show that $f:M\rightarrow P$ is a morphism of $G$-modules. To check this, recall (see beginning of Section \ref{S3}) that we are assuming that $\Spec R$ is connected and $R$ is a reduced finitely generated $\mathbb Z$-algebra. As each smooth $R$-algebra $S$ is still a reduced finitely generated $\mathbb Z$-algebra, by replacing $R$ with smooth $R$-algebras whose spectra are connected, it suffices to show that for every $g\in G(R)$ we have (cf. Equation (\ref{EQ01}))
\begin{equation}\label{EQ07}
f\circ\rho_M(R)(g)=\rho_{P}(R)(g)\circ f:M\rightarrow P.
\end{equation}
To check this, by giving up on the second recalled assumption on $R$, we can assume that $R=K=\overline{K}$ is an algebraically closed field and we will only use $K$. 

Let $\pmb{A}$ be the subgroup of $\pmb{\Aut}_M\times_{\Spec K} \pmb{\Aut}_{P}$ defined by the identity
$$\pmb{A}(K)=\{(g_1,g_2)\in \pmb{\Aut}_M(K)\times \pmb{\Aut}_{P}(K)|f\circ g_1=g_2\circ f\}.$$

Let $I$ be the image of the homomorphism 
$$(\rho_M,\rho_{P}):G\rightarrow \pmb{\Aut}_M\times_{\Spec K} \pmb{\Aut}_{P}.$$ Considering the short exact sequence $1\rightarrow\Ker(G\rightarrow I)\rightarrow G\rightarrow I\rightarrow 1$, as $\Ker(G\rightarrow I)$ is smooth over $\Spec K$ (due to Cartier's theorem), we get a short exact sequence of Lie algebras
$$0\rightarrow\Lie(\Ker(G\rightarrow I))\rightarrow\mathfrak g\rightarrow\Lie(I)\rightarrow 0.$$ As $f:L\rightarrow J$ is a morphism of $\mathfrak g$-modules and as $\mathfrak g\rightarrow\Lie(I)$ is surjective, we get that we have an inclusion
$$\Lie(I)\subset\Lie(\pmb{A}).$$
From this and \cite{Bor2}, Ch. II, Sect. 7, Subsect. 7.1 we get that $I$ is a subgroup of $\pmb{A}$ which implies that Equation (\ref{EQ07}) holds. Thus the functor (\ref{EQ02}) is surjective on morphisms. We conclude that Theorem \ref{T1} holds.\endproof

\section{Proof of Theorem \ref{T3}}\label{S4}

We will first prove the following basic lemma:

\begin{lemma}\label{L1} Let $I$ be a subgroup of a semisimple group $H$ of adjoint type over the spectrum of an algebraically closed field $\kappa$ such that $\dim(H/I)=1$. Then there exists an isomorphism $H\simeq\pmb{\PGL}_{2,\kappa}\times_{\Spec\kappa} H'$ which induces via restriction an isomorphism $I\simeq B_{2,\kappa}\times_{\Spec\kappa} H'$, where $B_{2,\kappa}$ is a Borel subgroup of $\pmb{\PGL}_{2,\kappa}$ and $H'$ is an arbitrary adjoint group over $\Spec\kappa$.\footnote{Another approach to prove this lemma due to Gabber is to use induction on the number of simple factors of $H$. One is reduced to the base of the induction case, so $H$ is simple, and it would suffice to prove that $\dim(H/I)$ is at least equal to the rank of $H$; such an inequality is well-known in characteristic $0$ but we could not find a reference for it in positive characteristic (however see \cite{LS}).}
\end{lemma}

\noindent
{\bf Proof:} We consider the connected smooth projective curve $C$ having $H/I$ as an open subscheme and let $g(C)$ be its genus. 

For simplicity, we define $\Aut(H/I)$ to be the reduced (smooth) subgroup of the group of automorphisms $\Aut(C)$ of $C$ that leave invariant the complement $C_0:=C\setminus (H/I)$. As for each field extension $\mu$ of $\kappa$, every automorphism of $(H/I)_{\mu}$ extends to an automorphism of $C_{\mu}$, the left multiplication action of $H(\mu)$ on $(H/I)_{\mu}$ induces an abstract homomorphism $H(\mu)\rightarrow\Aut(H/I)(\mu)$. Taking $\mu$ to be the field of fractions of $H$, we obtain a natural rational morphism from $H$ to $\Aut(H/I)$ and using translates it follows that it is defined everywhere. The reduced kernel of the resulting homomorphism $H\rightarrow\Aut(H/I)$ has a connected component $H'$ of the identity element which is the largest connected normal smooth subgroup of $H$ contained in $I$. Thus, as $H/H'$ is semisimple, we have inequalities 
\begin{equation}\label{EQ08}
3\le \dim(H/H')\le\dim(\Aut(H/I)).
\end{equation}

As $\Aut(H/I)\subset\Aut(C)$ and as the connected component $\Aut(C)^0$ of the identity element of $\Aut(C)$ is trivial if $g(C)\ge 2$, is an elliptic curve (thus abelian) if $g(C)=1$, and it is $\pmb{\PGL}_{2,\kappa}$ if $g(C)=0$, we conclude that $g(C)=0$ and we have a finite homomorphism $H/H'\rightarrow \pmb{\PGL}_{2,\kappa}$. As $H/H'$ is semisimple, by reasons of dimensions or by the classification of adjoint groups over $\kappa$, we get that $H/H'$ is isomorphic to either $\pmb{\PGL}_{2,\kappa}$ or $\pmb{\SL}_{2,\kappa}$. As $H$ is adjoint, the short exact sequence $1\rightarrow H'\rightarrow H\rightarrow H/H'\rightarrow 1$ splits. Thus $H\simeq H'\times_{\Spec\kappa} H/H'$ and we conclude that $H/H'\simeq\pmb{\PGL}_{2,\kappa}$. 

If $H/I$ is an affine rational curve, then $C\simeq \mathbb P^1_{\kappa}$ and the connected component $\Aut(H/I)^0$ of the identity element of $\Aut(H/I)$ is the subgroup of $\Aut(\mathbb P^1_{\kappa})\simeq\pmb{\PGL}_{2,\kappa}$ that fixes the finite nonempty set $C_0$; it follows that $\dim(\Aut(H/I))\le\dim(\pmb{\PGL}_{2,\kappa})-1=2$ which contradicts Inequality (\ref{EQ08}). 

Thus $H/I$ is projective isomorphic to $\mathbb P^1_{\kappa}$ which implies that $I/H'$ is a parabolic subgroup of $H/H'$, hence a Borel subgroup of $H/H'\simeq\pmb{\PGL}_{2,\kappa}$. The lemma follows from the last sentence and the isomorphisms $H/H'\simeq \pmb{\PGL}_{2,\kappa}$ and $H\simeq H'\times_{\Spec\kappa} H/H'$.\endproof

\medskip
To prove Theorem \ref{T3}, for an $H$-module $\diamond$, let $\ell_H(\diamond)$ be its length and let $\ell_{\Lie(H)}(\diamond)$ be its length as a $\Lie(H)$-module. We have a general inequality
\begin{equation}\label{EQ09}
\ell_H(\diamond)\le \ell_{\Lie(H)}(\diamond).
\end{equation}
As $\ell_H$ and $\ell_{\Lie(H)}$ are additive, from Inequality (\ref{EQ09}) we get immediately:

\begin{fact}\label{F1} 
If the Inequality (\ref{EQ09}) is an equality for $\diamond$, then it is an equality for each $H$-submodule or quotient of $\diamond$.
\end{fact}

We will use Lemma \ref{L1} to prove Theorem \ref{T3}, i.e., that $P$ is a semisimple $H$-module, by induction on $\ell:=\ell_H(P)=\ell_{\Lie(H)}(P)\in\mathbb Z_{\ge 0}$. The base of the induction for $\ell\in\{0,1\}$ is trivial. For $\ell\ge 2$ the passage from $\ell-1$ to $\ell$ goes as follows. Let $Q$ be a simple $H$-submodule of $P$: it is a simple $\Lie(H)$-module (by Fact \ref{F1}) and the $\Lie(H)$-module $P/Q$ is semisimple of the same length $\ell-1$ as the $H$-module $P/Q$. By the induction assumption, the $H$-module $P/Q$ is semisimple. Thus, to prove that the short exact sequence 
$$0\rightarrow Q\rightarrow P\rightarrow P/Q\rightarrow 0$$
splits, we can assume that $\ell=2$. 

As $\ell=2$, the $\Lie(H)$-module $P/Q$ is simple and we consider a simple $\Lie(H)$-submodule $N$ of $P$ which maps isomorphically onto $P/Q$. We have a direct sum decomposition $P=Q\oplus N$ of $\Lie(H)$-modules. 

We consider two cases as follows.

\smallskip\noindent
{\bf Case 1: the $\Lie(H)$-modules $N$ and $Q$ are not isomorphic.} Thus $P$ has only two simple $\Lie(H)$-submodules: $Q$ and $N$. As for all $h\in H(\kappa)$, we have $h\Lie(H)h^{-1}=\Lie(H)$, we get that $h(Q)$ and $h(N)$ are simple $\Lie(H)$-modules
for all $h\in H(\kappa)$. As $H$ is connected, from the last two sentences we get that for all $h\in H(\kappa)$ we have $Q=h(Q)$ and $N=h(N)$. This implies that both $Q$ and $N$ are $H$-submodules of $P$ which, as $\ell=2$, are simple. Thus $P=Q\oplus N$ is a semisimple $H$-module in this case. 

\smallskip\noindent
{\bf Case 2: the $\Lie(H)$-modules $N$ and $Q$ are isomorphic.} We fix a $\Lie(H)$-isomorphism $a:Q\rightarrow N$: it is unique up to multiplication by nonzero elements of $\kappa$. All simple $\Lie(H)$-submodules of $P$ are of the form 
$$Q_{[t_0:t_1]}:=\{t_0x+t_1a(x)|x\in Q\}\subset P=Q\oplus N$$ 
for a uniquely determined point $[t_0:t_1]\in\mathbb P^1_{\kappa}(\kappa)$. For instance, $Q=Q_{[1:0]}$ and $N=Q_{[0:1]}$. Similar to Case 1, for each field extension $\mu$ of $\kappa$ and for every $h\in H(\mu)$, $h(\mu\otimes_{\kappa} Q)$ is a simple $\Lie(H_{\mu})$-module and hence there exists a unique point $\delta(h)=[v_0:v_1]\in\mathbb P^1_{\kappa}(\mu)$ such that $$h(Q)=Q_{\delta(h)}:=\{v_0x+v_1a(x)|x\in\mu\otimes_{\kappa}Q\}.$$ 
An argument similar to the one involving $\mu$s in the proof of Lemma \ref{L1} shows that the association $h\rightarrow \delta(h)$ defines a morphism $\delta:H\rightarrow\mathbb P^1_{\kappa}$.

We show that the assumption that $\delta$ is nonconstant leads to a contradiction. For the stabilizer $I$ of $Q$ in $H$ we have $\dim(H/I)=\dim(\Im(\delta))=1$ and from Lemma \ref{L1} applied to the adjoint group $H^{\ad}$ of $H$ we get that there exists an isomorphism $H^{\ad}\simeq \pmb{\PGL}_{2,\kappa}\times_{\Spec\kappa} H'$ which induces via restriction an isomorphism $\Im(I\to H^{\ad})\simeq B_2\times_{\Spec\kappa} H'$, where $B_2$ is a Borel subgroup of $\pmb{\PGL}_{2,\kappa}$ and where $H'$ is an adjoint group over $\kappa$. This implies that $\delta$ is surjective and thus the simple $\Lie(H)$-submodules of $P$ are permuted transitively under the natural left multiplication action by $H(\kappa)$. But a simple $H$-submodule of $P$ is among the simple $\Lie(H)$-submodules of $P$ (by Fact \ref{F1}) and it is fixed by $H(\kappa)$, hence we reached a contradiction. 

Thus $\delta$ is constant of constant value $[1:0]$. Hence for all $h\in H(\kappa)$ we have $h(Q)=Q$ which implies that $Q$ is an $H$-submodule of $P$. The same applies to $N$. Thus $P=Q\oplus N$ is a semisimple $H$-module even in Case 2.

This ends the induction and the proof of Theorem \ref{T3}.\endproof

\section{Proof of Theorem \ref{T2}}\label{S5}

In this section we assume that $R$ is a faithfully flat $\mathbb Z_{(p)}$-algebra and a normal integral domain. 

Theorem \ref{T2} (a) is proved in Subsection \ref{S51}. Theorem \ref{T2} (b) is proved in Subsection \ref{S52}. If $G$ is split, then the hypotheses of Theorem \ref{T2} (b) hold: as $D$ we can take $\mathbb Z_{(p)}$, see the existence of $G_{\mathbb Z}$ and $V_{\mathbb Q}$ in Subsection \ref{S22}. Thus, Theorem \ref{T2} (c) follows directly from Theorems \ref{T2} (a) and (b). 

We recall the following well-known fact.

\begin{fact}\label{F2}
Let $\mathcal S$ be an affine smooth scheme over the spectrum of a discrete valuation ring $\mathcal D$ with uniformizer $\varpi$ and residue field $k$. Let $\mathcal M$ be a free $\mathcal D$-module which is a $\mathcal S$-module (i.e., it is equipped with a homomorphism $\mathcal S\rightarrow\pmb{\Aut}_{\mathcal M}$). Let $\mathcal L$ be a $\mathcal D$-submodule of $\varpi^{-1}\mathcal M$ which contains $\mathcal M$. Then $\mathcal L$ is a $\mathcal S$-module (resp. a $\Lie(\mathcal S)$-module) if and only if $\mathcal L/\mathcal M$ is a $\mathcal S_k$-submodule (resp. is a $\Lie(\mathcal S_k)$-submodule) of $\varpi^{-1}\mathcal M/\mathcal M$.
\end{fact} 

\noindent
{\bf Proof:} The `only if' parts and the case of Lie algebras are obvious, hence it suffices to check that if $\mathcal L/\mathcal M$ is a $\mathcal S_k$-module, then $\mathcal L$ is a $\mathcal S$-module. Writing $\mathcal S=\Spec\mathcal A$, this is equivalent to checking that the comultiplication $\mathcal D$-linear map $\nabla:\mathcal M\rightarrow \mathcal A\otimes_{\mathcal D} \mathcal M$ is such that $\nabla(\mathcal L)\subset\mathcal A\otimes_{\mathcal D} \mathcal L$. To check this we can assume that $k$ is algebraically closed (it suffices to be infinite) and $\mathcal D$ is complete and we will check directly (i.e., without mentioning $\nabla$ again) that $\mathcal L$ is a $\mathcal S$-module. As $\mathcal L/\mathcal M$ is a $\mathcal A_k$-submodule of $\varpi^{-1}\mathcal M/\mathcal M$, we get that for each $h\in\mathcal S(\mathcal D)$ we have $h(\mathcal L)=\mathcal L$. From the last two sentences and \cite{V1}, Prop. 3.1.2.1 a) we get that the homomorphism $\mathcal S\rightarrow\pmb{\Aut}_{\mathcal M}$ over the spectrum of the field of fractions of $\mathcal D$ extends to a homomorphism $\mathcal S\rightarrow\pmb{\Aut}_{\mathcal L}$, thus $\mathcal L$ is a $\mathcal S$-module.\endproof

\subsection{Proof of Theorem \ref{T2} (a)}\label{S51}

Let $Z:=\Ker(G^{\sc}\to G)$; it is a finite flat group scheme over $\Spec R$ of multiplicative type which is contained in the center of $G^{\sc}$ and which, if $G$ is split, is the kernel of the induced homomorphisms between split maximal tori. In particular, the Zariski (or the schematic) closure of $Z_K$ in $G^{\sc}$ is $Z$ itself. Thus, if $M$ is a $G^{\sc}$-module such that we can identify $K\otimes_R M=V$ as $G^{\sc}_K$-modules, then the kernel of the homomorphism $\rho_M:G^{\sc}\to\pmb{\Aut}_M$ contains $Z_K$ and therefore it contains $Z$; hence $M$ is in fact a $G$-module. Thus we can identify $\Lat_{G^{\sc}}(V)=\Lat_G(V)$. Based on this and the inclusions $\Lat_G(V)\subset\Lat_{\mathfrak g}(V)\subset \Lat_{\Lie(G^{\sc})}(V)\supset \Lat_{G^{\sc}}(V)$, it suffices to prove that Theorems \ref{T2} (a) and (b) hold in the case when $G=G^{\sc}$ is simply connected. 

To prove that Theorem \ref{T2} (a) holds, as in the beginning of Section \ref{S3}, using inductive limits and working in the \'etale topology of $\Spec R$, we can assume that $R$ is also noetherian and that $G$ is split. Based on Case 2 of Subsection \ref{S31}, as in Case 3 of Subsection \ref{S31} we argue that Theorem \ref{T2} (a) holds provided it holds for discrete valuation rings of mixed characteristic $(0,p)$. Thus we can also assume that $R=D$ is a discrete valuation ring of mixed characteristic $(0,p)$; let $\varpi$ be a uniformizer of it. Let $k:=D/(\varpi)$: it is a field of characteristic $p$. Let $H:=G_k$ and let $\mathfrak h:=\mathfrak g/\varpi\mathfrak g=k\otimes_D \mathfrak g$.

Let $L\in\Lat_{\mathfrak g}(V)$; we have $V=K\otimes_D L$. With the notation of Subsections \ref{S22} and \ref{S23} for $G=G^{\sc}$, if we have $L\in\Lat_{G_i}(V)$ for all $i\in\{1,\ldots,n\}$, then we obtain homomorphisms $\rho_{L,i}:G_i\rightarrow\pmb{\Aut}_L$ whose fibers over $K$ are restrictions of $\rho_V:G_K\rightarrow \pmb{\Aut}_V$. These homomorphisms over $\Spec D$ commute as they commute over $K$. This implies that their product defines a homomorphism $\rho_L:G\rightarrow \pmb{\Aut}_L$ which extends $\rho_V$ and hence we have $L\in\Lat_G(V)$. Thus to prove Theorem \ref{T1} (a) we can assume that $n=1$, i.e., the split adjoint group scheme $G^{\ad}$ of $G=G^{\sc}=G_1$ has simple geometric fibers. 

As in Case 2 of Subsection \ref{S31}, based on \cite{J}, Part I, Ch. X, Lem. of Sect. 10.4 we get the existence of a lattice $M$ of $V$ with respect to $R$ which is a $G$-module. So $M/\varpi M$ is an $H$-module. 

To prove that $L$ is a $G$-module, we can assume that $k$ is algebraically closed and we can replace $L$  by $\varpi^{\tau}L$ with $\tau\in\mathbb Z$. Thus we can assume that $M\subset L$ but $\varpi^{-1}M\not\subset L$. Let $s\ge 0$ be the smallest integer such that $L\subset \varpi^{-s}M$. We will prove by induction on $s\in\mathbb Z_{\ge 0}$ that, regardless of what the $G$-module $M$ is, $L$ is a $G$-module, i.e., $L\in\Lat_G(V)$. The base of the induction is trivial: if $s=0$, then $L=M$ is a $G$-module. 

For $s\in\mathbb N$ the passage from at most $s-1$ to $s$ goes as follows. If $s\ge 2$, then $M\subset \varpi L+M\subset\varpi^{-s+1}M$ are inclusions between $\mathfrak g$-modules and hence by induction applied first with $(L,s)$ replaced by $(\varpi L+M,s-1)$ we get that $\varpi L+M$ is a $G$-module and applied second with $(M,s)$ replaced by $(\varpi L+M,1)$ we get that $L$ is a $G$-module. Thus we can assume that $s=1$, i.e., we have $M\subsetneq L\subsetneq \varpi^{-1}M$. Let 
$$\mathfrak n:=L/M\subset \mathfrak m:=\varpi^{-1}M/M\simeq M/\varpi M;$$ 
it is a nonzero $\mathfrak h$-module. 

We will prove using a second induction on the length $t\in\mathbb N$ of the $\mathfrak h$-module $\mathfrak n$ that, regardless of what the $G$-module $M$ is, $L$ is a $G$-module. Let $\mathfrak p$ be a $\mathfrak h$-submodule of $\mathfrak n$ of length $t-1$: thus the $\mathfrak h$-module $\mathfrak n/\mathfrak p$ is simple. Let $M_+$ be the inverse image of $\mathfrak p$ via the $D$-linear map $\varpi^{-1}M\rightarrow\varpi^{-1}M/M=\mathfrak m$. If $t=1$, then $M_+=M$ is a $G$-module. Thus, if $t\ge 2$ and the statement is true for $\le t-1$, then by the (second) induction assumption we get first that $M_+$ is a $G$-module and second, by replacing $(M,t)$ by $(M_+,1)$, that $L$ itself is a $G$-module. Hence to end the proof of both inductions we can assume that not only $s=1$ but we also have $t=1$. So $\mathfrak n$ is a simple $\mathfrak h$-module. 

For a maximal torus $T$ of $G$ which is split, the weights of the action of $T$ on $M$ are the same as the weights of the action of $T_k$ on $\mathfrak m$. Based on this, as statement \circled{\textup{2}} holds, for each composition series of the $H$-module $\mathfrak m$, the simple factors are irreducible $H$-modules associated to highest weights $\sum_{l=1}^{r_1} \gamma_{1,l}\omega_{1,l}$ with the property that for all $l\in\{1,\ldots,r_1\}$ we have $\gamma_{1,l}\in\{0,\ldots,p-1\}$ and hence are simple $\mathfrak h$-modules (see \cite{Bor1}, Thm. 6.4). This implies that the $H$-module $\mathfrak m$ and the $\mathfrak h$-module $\mathfrak m$ have the same length, and, by Fact \ref{F1}, the same holds for each $H$-submodule $\mathfrak p$ of $\mathfrak m$. 

We take $\mathfrak p$ to be the $H$-submodule of $\mathfrak m$ generated by $\mathfrak n$. As $k$ is algebraically closed, we have an identity
$$\mathfrak p=\sum_{h\in H(k)} h(\mathfrak n)$$
of $k$-vector spaces. As $h$ normalizes $\mathfrak h$, each $h(\mathfrak n)$ is a simple $\mathfrak h$-module and therefore $\mathfrak p$, being a sum of simple $\mathfrak h$-modules, is a semisimple $\mathfrak h$-submodule of $\mathfrak m$. From Theorem \ref{T3} we get that $\mathfrak p$ is a semisimple $H$-module. 

Writing $\mathfrak p=\oplus_{u=1}^b \mathfrak p_u$ as a direct sum of simple $H$-modules, from Fact \ref{F1} we get that each $\mathfrak p_u$ is a simple $\mathfrak h$-module. From this and the fact that $\mathfrak p$ is the $H$-submodule of $\mathfrak m$ generated by $\mathfrak n$, we get that the $\mathfrak h$-module $\mathfrak n$ projects isomorphically onto each $\mathfrak p_u$. Hence the $\mathfrak h$-module $\mathfrak p$ is isomorphic to $b\mathfrak n:=\oplus_{u=1}^b \mathfrak n$. If $\omega_u$ is the highest weight of the $H$-module $\mathfrak p_u$, then as the isomorphism class of the $\mathfrak h$-module $\mathfrak p_u$ does not depend on $u$, we easily get that $\omega_u\in\Omega_{1}:=\{\sum_{l=1}^{r_1} \gamma_{1,l}\omega_{1,l}|\gamma_{1,1},\ldots,\gamma_{1,r_1}\in\{0,1,\ldots,p-1\}\},$ does not depend on $u$ (this is also proved in \cite{Bor1}, Subsect. 6.6). Thus we have $\mathfrak p=b\mathfrak p_1$ as $H$-modules as well as  $\mathfrak h$-modules and this implies that a $k$-vector subspace of $\mathfrak p$ is an $H$-module if and only if it is a  $\mathfrak h$-module. Therefore $\mathfrak n$ is an $H$-module (and in particular we have $\mathfrak p=\mathfrak n$ and $b=1$). This implies that $L$ is a $G$-module (see Fact \ref{F2}). This ends the proof of both inductions and hence of Theorem \ref{T2} (a).\endproof

\subsection{Proof of Theorem \ref{T2} (b)}\label{S52}

Using the contrapositive, it suffices to show that if statement \circled{\textup{2}} does not hold, then there exists $L\in\Lat_{\mathfrak g} (V)$ which is not a $G$-module. Considering pullbacks via $\Spec R\rightarrow\Spec D$ (i.e., the tensorization of elements of $\Lat_{\Lie(G_D)}(V_{\Frac{D}})$ over $D$ with $R$), to find such an $L$, we can assume that $R=D$ is a discrete valuation ring of mixed characteristic $(0,p)$. Let $\varpi$, $k$ and $M$ be as in Subsection \ref{S51}: we will only use $G$ (as $G=G^{\sc}=G_D$), $D$ (as $R=D$), and $V$ (as $V=V_K$). 

The factors of a composition series of the $H:=G_k$-module $\mathfrak m\simeq M/\varpi M$ do not depend on the choice of the $G$-module $M$, see \cite{J}, Part I, Ch. X, Sects. 10.7 and 10.9. Let $D^{\h}$ be the henselization of $D$. 

We consider two disjoint cases as follows.

\smallskip\noindent
{\bf Case 1: $H$ is split.} As $H$ is split, the affine smooth scheme $\mathcal T_G$ over $\Spec D$ that parametrizes maximal tori of $G$ (see \cite{SGA3-2}, Exp. XII, Cors. 1.10 and 5.4) has a $k$-valued point defining a split maximal torus of $H$ and therefore, due to the smoothness of $\mathcal T_G$, it lifts to a $D^{\h}$-valued point of $\mathcal T_G$. Thus $G_{D^{\h}}:=\Spec D^{\h}\times_{\Spec D} G$ has a maximal torus whose fiber over $\Spec k$ is split and hence, as $D^{\h}$ is henselian, it is split. Thus $G_{D^{\h}}$ is split. This implies that there exists a $D$-subalgebra $D'$ of $\overline{K}$ which is \'etale and a discrete valuation ring of residue field $k$ and which is such that $G_{D'}$ is split. 

Let $K':=\Frac(D')$. If $T'$ is a maximal torus of $G_{D'}$ which is split, then the weights of the action of $T'_{K'}$ on $K'\otimes_K V=K'\otimes_D M$ and of the action of $T'_k$ on $\mathfrak m$ are the same. As statement \circled{\textup{2}} does not hold, we deduce that the composition series of the $H=\prod_{i=1}^n G_{i,k}$-module $\mathfrak m$ has a simple factor $\mathfrak n$ which, up to isomorphism, is a tensor product $\otimes_{i=1}^n \mathfrak n_i$, where each $\mathfrak n_i$ is a simple $G_{i,k}$-module of highest weight $w_i$, and there exists $i_0\in\{1,\ldots,n\}$ such that we can write
$$w_{i_0}=\sum_{l=1}^{r_{i_0}} c_{i_0,l}\omega_{i_0,l}$$
with all $c_{i_0,l}\in\mathbb Z_{\ge 0}$ but there exists $l_0\in\{1,\ldots,r_{i_0}\}$ such that $c_{i_0,l_0}\ge p$. This implies that we can write
$$w_{i_0}=\sum_{t=0}^q p^tw_{i_0,t}$$
with $q\in\mathbb N$, and $w_{i_0,q}\neq 0$, and
$$w_{i_0,0},\ldots,w_{i_0,q}\in\Omega_{i_0}:=\{\sum_{l=1}^{r_{i_0}} \gamma_{i_0,l}\omega_{i_0,l}|\gamma_{i_0,1},\ldots,\gamma_{i_0,r_{i_0}}\in\{0,1,\ldots,p-1\}\}.$$ The key point is (see \cite{St}, Thm. 1.1; see also \cite{Bor1}, Thm. 7.5 (i)) that we have a tensor product decomposition
$$\mathfrak n_{i_0}\simeq\otimes_{t=0}^q \mathfrak n_{i_0,t}^{(p^t)}$$
to be viewed as an identification, where $\mathfrak n_{i_0,t}$ is the simple $G_{i_0,k}$-module associated to the highest weight $w_{i_0,t}$ and we view the $G_{i_0,k}^{(p^t)}$-module $\mathfrak n_{i_0,t}^{(p^t)}$ as a $G_{i_0,k}$-module via the functorial Frobenius homomorphism $G_{i_0,k}\rightarrow G_{i_0,k}^{(p^t)}$. We recall that $G_{i_0,k}^{(p^t)}$ is the pullback of $G_{i_0,k}$ via the morphism $\Spec k\rightarrow\Spec k$ defined by the Frobenius endomorphism $\Fr_t:k\rightarrow k$  that maps $x$ to $x^{p^t}$ and that $\mathfrak n_{i_0,t}^{(p^t)}:= k\otimes_{\Fr_t,k} \mathfrak n_{i_0,t}$. For $t\in\{1,\ldots,q\}$, $\Lie(G_{i_0,k})$ acts trivially on $\mathfrak n_{i_0,t}^{(p^t)}$. 

We consider two $H$-submodules $\mathfrak p_0\subset \mathfrak p_1$ of $\mathfrak m$ such that as $H$-module $\mathfrak p_1/\mathfrak p_0$ is (isomorphic to) such an $H$-module $\mathfrak n$ and the $H$-module $\mathfrak p_0$ has the smallest length. By replacing the $G$-module $M$ by the inverse image of $\mathfrak p_0$ (see Fact \ref{F2}) via the composite $D$-linear map
$$\nu:M\rightarrow M/\varpi M\simeq\mathfrak m,$$ we can assume that $\mathfrak p_0=0$. Thus $\mathfrak p_1=\mathfrak p_1/\mathfrak p_0$ is a simple $H$-module. Therefore we can assume that $\mathfrak n=\mathfrak p_1$ is a simple $H$-submodule of $\mathfrak m$.

As $w_{i_0,q}\neq 0$, we have $\dim_k(\mathfrak n_{i_0,q})\ge 2$ and thus there exists a nonzero proper $k$-vector subspace $Q_{i_0,q}$ of $\mathfrak n_{i_0,q}^{(p^q)}$: it is a trivial $\Lie(G_{i_0,k})$-module which is not a $G_{i_0,k}$-module. The $k$-vector subspace 
$$\mathfrak q_{i_0}:=\mathfrak n_{i_0,0}^{(p^0)}\otimes_k \mathfrak n_{i_0,1}^{(p^1)}\otimes_k\cdots\otimes_k\mathfrak n_{i_0,q-1}^{(p^{q-1})}\otimes_k Q_{i_0,q}$$ 
of $\mathfrak n_{i_0}$ is a $\Lie(G_{i_0,k})$-module which is not a $G_{i_0,k}$-module. Defining 
$$\mathfrak q:=\mathfrak n_1\otimes_k\cdots\otimes_k\mathfrak n_{i_0-1}\otimes_k \mathfrak q_{i_0}\otimes_k\mathfrak n_{i_0+1}\otimes_k\cdots\otimes_k\mathfrak n_n,$$ 
we get that $\nu^{-1}(\mathfrak q)\in\Lat_{\mathfrak g}(V)$ but (see Fact \ref{F2}) $\nu^{-1}(\mathfrak q)\notin\Lat_G(V)$ as $\mathfrak q$ is not an $H$-module.

\smallskip\noindent
{\bf Case 2: $H$ is not split.} Let $k'$ be a finite separable field extension of $k$ such that $H_{k'}:=\Spec k'\times_{\Spec k} H$ is split. We can assume that the field extension $k\to k'$ is Galois. The Galois group $\Gamma:=\Gal(k'/k)$ acts naturally on the set $\{1,\ldots,n\}$ that indexes the absolutely simple factors of $\Spec k'\times_{\Spec k} G_k^{\ad}=G^{\ad}_{k'}=\prod_{i=1}^n G^{\ad}_{i,k'}$. 

Let $D'$ be a discrete valuation ring which is a finite flat $D$-algebra with the property that we have an identity $D'/\varpi D'=k'$. Let $K':=\Frac(D')$; we have $K'=D'\otimes_D K=D'[1/\varpi]$. Let $\mathfrak g':=\Lie(G_{D'})=D'\otimes_D \mathfrak g$. 

As $\Gamma$ is canonically identified with $\Aut(K'/K)$ and $\Aut(D'/D)$, it acts naturally on $K'\otimes_K V$ and $D'\otimes_D M$.

Based on Case 1 applied to $G_{D'}$ and the $G_{K'}$-module $K'\otimes_K V$, we get the existence of a lattice $L'\in\Lat_{\mathfrak g'}(K'\otimes_K V)$ which is not a $G_{D'}$-module. Let $L:=L'\cap V$, the intersection being taken inside $K'\otimes_K V$. As $L'$ is a $\mathfrak g'$-module, it is also a $\mathfrak g$-module and we conclude that $L\in\Lat_{\mathfrak g}(V_K)$. If $L$ is a $G$-module, then $D'\otimes_D L$ is a $G_{D'}$-module and hence we have $L'\neq D'\otimes_D L$. Thus to end the proof in this case it suffices to show that we can choose $L'$ such that we have $L'=D'\otimes_D L$. 

Based on Case 1 applied over $D'$ to the $G_{D'}$-module $D'\otimes_D M$, we consider $H_{k'}$-submodules $\mathfrak p_0'\subset\mathfrak p_1'\subset (D'\otimes_D M)/(\varpi D'\otimes_D M)=k'\otimes_k \mathfrak m$ such that $\mathfrak n':=\mathfrak p'_1/\mathfrak p'_0$ is a simple $H_{k'}$-module that has the same property as $\mathfrak n$ of Case 1. We can assume that $\mathfrak p'_0$ is such that its length as an $H_{k'}$-module is the smallest. Due to this, by considering $\mathfrak r^{\prime}_0:=\sum_{\gamma\in \Gamma} \gamma(\mathfrak p_0)$, which, due to Galois descent, is of the form $k'\otimes_k \mathfrak r_0$ with $\mathfrak r_0$ an $H$-submodule of $\mathfrak m$, and $\mathfrak r'_1:=\mathfrak p_1+\mathfrak r^{\prime}_0$, we have an $H_{k'}$-isomorphism $\mathfrak n'_0\simeq \mathfrak r'_1/\mathfrak r'_0$. Thus, as in Case 1, by replacing $M$ with $\nu^{-1}(\mathfrak r_0)$, we can assume that $\mathfrak p'_0=0$ and hence we have a simple $H_{k'}$-module $\mathfrak n'$ of $k'\otimes_k \mathfrak m$ which is the analogue of $\mathfrak n$.  

As in Case 1 we get a tensor product decomposition $\mathfrak n'=\prod_{i=1}^n \mathfrak n'_i$, an element $i_0\in\{1,\ldots,n\}$, and a second product decomposition $\mathfrak n'_{i_0}\simeq\otimes_{t=0}^q (\mathfrak n'_{i_0,t})^{(p^t)}$. The stabilizer $\Gamma_{i_0}$ of $i_0$ in $\Gamma$ acts naturally on the set of dominant weights of $G_{i_0,k'}$ and as such let $\Gamma_{i_0}^-$ be the subgroup of $\Gamma_{i_0}$ that fixes $w_{i_0,q}$. Thus the $G_{i_0,k'}$-module $\mathfrak n'_{i_0,q}$ is defined over the subfield of $k'$ fixed by $\Gamma_{i_0}^-$ and as such there exists a nonzero proper $k'$-vector subspace $Q'_{i_0,q}$ of $(\mathfrak n'_{i_0,q})^{(p^q)}$ left invariant by $\Gamma_{i_0}^-$. We use it to define an element 
$$M':=(1_{D'}\otimes_D \nu)^{-1}(\mathfrak q')\in\Lat_{\mathfrak g'}(K'\otimes_K V)$$ 
which (see Fact \ref{F2}) is not a $G_{D'}$-module. Then we can take 
$$L':=\sum_{\gamma\in\Gamma} \gamma(M').$$ 
As $L'$ is $\Gamma$-invariant, using Galois descent we get that we have $L'= D'\otimes_D L$. As $\mathfrak g'$ is $\Gamma$-invariant, each $\gamma(M')$ is a $\mathfrak g'$-module. As $\Gamma_{i_0}^{-}$ leaves $M'$ invariant, we have $L'=\sum_{\gamma\in\Gamma/\Gamma_{i_0}^{-}} \gamma(M')$ and it follows that $L'$ is not a $G_{D'}$-module. 

We conclude that Theorem \ref{T2} (b) holds.

\section{Proof of Theorem \ref{T4}}\label{S6}

Due to the perfect assumptions we have a direct sum decomposition
$$\End_R(\mathcal M)=\Lie(\mathcal H_K)\cap\End_R(\mathcal M)\oplus [\Lie(\mathcal H_K)\cap\End_R(\mathcal M)]^{\perp},$$
where $[\Lie(\mathcal H_K)\cap\End_R(\mathcal M)]^{\perp}$ is the perpendicular of $\Lie(\mathcal H_K)\cap\End_R(\mathcal M)$ with respect to $\mathcal B$. Thus the $R$-module underlying the Lie algebra $\Lie(\mathcal H_K)\cap\mathfrak g\mathfrak l_{\mathcal M}$ is projective. Moreover, the Killing form of $\Lie(\mathcal H_K)\cap\mathfrak g\mathfrak l_{\mathcal M}$ is a perfect bilinear map as it times a unit of $R$ is so. As $\char(K)=0$, we have $\Lie(\mathcal H_K)=\Lie(\mathcal H_K^{\ad})$. Based on the last two sentences, from \cite{V3},  Cor. 1.3 we get that there exists a unique adjoint group scheme $\mathcal H^{\ad}$ over $\Spec R$ which extends $\mathcal H_K^{\ad}$ and whose Lie algebra is $\Lie(\mathcal H_K)\cap\mathfrak g\mathfrak l_{\mathcal M}$. 

Let $\mathcal G$ be the simply connected semisimple group scheme cover of $\mathcal H^{\ad}$; we have $\mathcal G^{\ad}=\mathcal H^{\ad}$ and our notation matches (i.e., the fiber of $\mathcal G$ over $\Spec K$ is the `initial' $\mathcal G_K$). Let $\mathcal H$ be the normalization of $\mathcal H^{\ad}$ in the field of fractions of $\mathcal H_K$. From \cite{V3}, Lem. 2.3.1 we get that $\mathcal H$ has a unique structure of a semisimple group scheme over $\Spec R$ which extends $\mathcal H_K$ and the morphism $\mathcal H\to\mathcal H^{\ad}$ is in fact a central isogeny. Clearly, $\mathcal G$ is also the simply connected semisimple group scheme cover of $\mathcal H$. 

From \cite{V3}, Prop. 3.5 (b) we get that the homomorphisms $\mathcal G\to\mathcal H$ and $\mathcal H\to\mathcal H^{\ad}$ are \'etale. Thus we have identifications 
$$\Lie(\mathcal G)=\Lie(\mathcal H)=\Lie(\mathcal H^{\ad})=\Lie(\mathcal H_K)\cap\mathfrak g\mathfrak l_{\mathcal M}$$ and the kernels $\mathcal K:=\Ker(\mathcal G\to\mathcal H)$ and $\Ker(\mathcal H\to\mathcal H^{\ad})$ are finite \'etale group schemes over $\Spec R$. 

As $\mathcal M$ is a $\Lie(\mathcal H_K)\cap\mathfrak g\mathfrak l_{\mathcal M}$-module and hence also a $\Lie(\mathcal G)$-module, from Theorem \ref{T2} (a) we get that $\mathcal M$ is a $\mathcal G$-module. The kernel of the resulting homomorphism $\mathcal G\to\pmb{\Aut}_{\mathcal M}$ contains $\mathcal K_K$ and hence contains $\mathcal K$. This implies that $\mathcal G\to\pmb{\Aut}_{\mathcal M}$ factors through a homomorphism $\rho:\mathcal H\to \pmb{\Aut}_{\mathcal M}$.

Using a limit argument as in the beginning of Section \ref{S3}, to check that $\rho$ is a closed embedding we can assume that the normal domain $R$ is also noetherian.

As $\rho_K$ is a closed embedding, from \cite{V2}, Thm. 1.1 (c) we get that $\rho$ is finite over the spectrum of each local ring of $R$ which is a discrete valuation ring. This implies that there exists an open subscheme $Y$ of $\Spec R$ which contains all points of $\Spec R$ of codimension in $\Spec R$ at most $1$ and such that $\rho_Y$ is finite. Based on this, as $p>2$, from \cite{V3}, Prop. 5.1 we get that $\rho$ itself is a closed embedding. This implies that the cokernel $\mathfrak g\mathfrak l_{\mathcal M}/Lie(\mathcal H)$  of the inclusion $\Lie(\mathcal H)\to\mathfrak g\mathfrak l_{\mathcal M}$ has constant rank over the points of $\Spec R$, so it is a finitely generated projective module over $R$. Thus $\Lie(\mathcal H)$ is a direct summand of $\mathfrak g\mathfrak l_{\mathcal M}$ and hence we have $\Lie(\mathcal H)=\Lie(\mathcal H_K)\cap\mathfrak g\mathfrak l_{\mathcal M}$. Therefore Theorem \ref{T4} holds.\endproof

\medskip\noindent
{\bf Acknowledgement.}
The authors confirm that the data supporting the findings of this study are available within the article. This paper has no data associated to it and no conflict of interest. We would like to thank Ofer Gabber for pointing out the essence of the footnote of Section \ref{S4}. The first author would like to thank SUNY at Binghamton for years of support during his gradate studies and the second author for introducing him to the topics of this paper, providing valuable advice and constant support, and always being willing to answer questions. The second author would like to thank SUNY at Binghamton for good working conditions.

\hbox{}
\hbox{Micah Loverro,\;\;\;E-mail: m.loverro@gmail.com}
\hbox{}

\bigskip

\hbox{Adrian Vasiu,\;\;\;E-mail: adrian@math.binghamton.edu}
\hbox{Address: Department of Mathematical Sciences, Binghamton University,}
\hbox{P. O. Box 6000, Binghamton, New York 13902-6000, U.S.A.}

\end{document}